\documentclass{amsart}
\usepackage{amsmath}
\usepackage{amssymb}
\usepackage{mathrsfs}
\usepackage{amscd}
\usepackage[all,knot,arc]{xy}\usepackage{graphicx}
\usepackage{subfigure}
\usepackage{verbatim}

\usepackage[plainpages]{hyperref}

\newtheorem{theorem}{Theorem}[section] 

\newtheorem{corollary}[theorem]{Corollary}

\newtheorem{proposition}[theorem]{Proposition} 

\theoremstyle{definition}

\newtheorem{remark}[theorem]{Remark}
\newtheorem*{ack}{Acknowledgements}

\numberwithin{equation}{section}

\numberwithin{figure}{section}

\newcounter{first}
{\end{list}}

\newif\ifpic
\picfalse    
\pictrue   

\newcommand{\A}{\ensuremath{\mathcal A}}
\newcommand{\MM}{{\ensuremath{\mathfrak M}}}
\newcommand{\C}{\ensuremath{\mathbb C}}
\newcommand{\RR}{\ensuremath{\mathcal R}}
\newcommand{\Z}{\ensuremath{\mathbb Z}}
\newcommand{\Q}{\ensuremath{\mathbb Q}}
\newcommand{\bS}{\ensuremath{\mathbb S}}
\newcommand{\tv}{\operatorname{tv}}
\newcommand{\Mod}{\operatorname{Mod}}
\newcommand{\Aut}{\operatorname{Aut}}
\newcommand{\Span}{\operatorname{span}}

\renewcommand{\k}{\ensuremath{\Bbbk}}
\renewcommand{\bar}{\overline}
   
\newcommand{\PSL}{\operatorname{PSL}}
\newcommand{\Aff}{\operatorname{Aff}}

\begin{document}
\title[Pure braid groups are not residually free]{Pure braid groups are not residually free}
\author[Daniel C. Cohen]{Daniel C. Cohen$^\dag$}
\address{Department of Mathematics, Louisiana State University, Baton Rouge, Louisiana 70803}
\email{\href{mailto:cohen@math.lsu.edu}{cohen@math.lsu.edu}}
\urladdr{\href{http://www.math.lsu.edu/~cohen/}
{www.math.lsu.edu/\char'176cohen}}
\thanks{{$^\dag$}Partially supported by Louisiana Board of Regents grant NSF(2010)-PFUND-171}

\author{Michael Falk}
\address{Department of Mathematics and Statistics, Northern Arizona University, Flagstaff, Arizona 86011-5717}
\email{\href{mailto:michael.falk@nau.edu}{michael.falk@nau.edu}}
\urladdr{\href{http://www.cefns.nau.edu/~falk/}
{www.cefns.nau.edu/\char'176falk}}
\author{Richard Randell}
\address{Department of Mathematics, University of Iowa, Iowa City, Iowa 52242}
\email{\href{mailto:randell@math.uiowa.edu}{randell@math.uiowa.edu}}
\urladdr{\href{http://www.math.uiowa.edu/~randell/}
{www.math.uiowa.edu/\char'176randell}}

\begin{abstract}
We show that the Artin pure braid group $P_n$ is not residually free for $n\ge 4$. 
Our results also show that the corank of $P_n$ is equal to $2$ for $n \ge 3$.
\end{abstract}

\maketitle

\section{Introduction}
A group $G$ is 
residually free if for every $x \neq 1$ in $G$, there is a homomorphism $f$ from $G$ to a free group $F$ so that $f(x)\neq 1$ in $F$. Equivalently, $G$ embeds in a product of free groups (of finite rank).  Examples of residually free groups include the fundamental groups of orientable surfaces.  In this note, we show that the Artin pure braid group, the kernel $P_n=\ker(B_n \to \Sigma_n)$ of the natural map from the braid group to the symmetric group, is not residually free for $n\ge 4$. 
(It is easy to see that the pure braid groups $P_2$ and $P_3$ are residually free.) 
We also classify all epimorphisms from the pure braid group to free groups, and determine the corank of the pure braid group. 
For $n\geq 5$, the fact that $P_n$ is not residually free was established independently by L.~Paris (unpublished), see Remark~\ref{rem:paris}. 

For $n\ge 3$, the braid groups themselves are not residually free.  Indeed, the only nontrivial two-generator residually free groups are $\Z$, $\Z^2$, and $F_2$, the nonabelian free group of rank two, see Wilton \cite{Wil08}.  
Since $B_n$ can be generated by two elements for $n\ge 3$, it is not residually free.  
(For $n=2$, $B_2=\Z$ is residually free.)

If $G$ is a group which is not residually free, then any group $\tilde{G}$ with a subgroup isomorphic to $G$ cannot be residually free.  Consequently, the (pure) braid groups are ``poison'' groups for residual freeness.  In particular, a group with a subgroup isomorphic to the $4$-strand pure braid group $P_4$ or the $3$-strand braid group $B_3$ is not residually free.  
Since $P_4<P_n$ for every $n\geq 4$, our main result follows from the special case $n=4$. 
Moreover, the same observation 
enables us to show that a number of other groups are not residually free.  These include (pure) braid groups of orientable surfaces, the (pure) braid groups associated to the full monomial groups, and a number of irreducible (pure) Artin groups of finite type.

This research was motivated by our work in \cite[\S 3]{CFR10}, which implies the residual freeness of fundamental groups of the complements of certain complex hyperplane arrangements. 
In particular, the proof of the assertion in the last sentence of \cite[Example 3.25]{CFR10} gives the last step in the proof of Theorem \ref{thm:notresfree} below. 

\section{Automorphisms of the pure braid group}

Let $B_n$ be the Artin braid group, with generators $\sigma_1,\dots,\sigma_{n-1}$ and relations $\sigma_i\sigma_{i+1}\sigma_i=\sigma_{i+1}\sigma_i\sigma_{i+1}$ for $1\le i\le n-2$, and $\sigma_{j}\sigma_i=\sigma_{i}\sigma_j$ for $|j-i|\ge 2$.  The Artin pure braid group $P_n$ has generators
\[
A_{i,j}=\sigma_{j-1}^{}\cdots\sigma_{i+1}^{}\sigma_i^{2}\sigma_{i+1}^{-1}\cdots\sigma_{j-1}^{-1}
=\sigma_{i}^{-1}\cdots\sigma_{j-2}^{-1}\sigma_{j-1}^{2}\sigma_{j-2}^{}\cdots\sigma_{i}^{},
\]
and relations
\begin{equation} \label{eq:purebraidrels}
A_{r,s}^{-1}A_{i,j}^{}A_{r,s}^{}=
\begin{cases}
A_{i,j}^{}&\text{if $i<r<s<j$,}\\
A_{i,j}^{}&\text{if $r<s<i<j$,}\\
A_{r,j}^{}A_{i,j}^{}A_{r,j}^{-1}&\text{if $r<s=i<j$,}\\
A_{r,j}^{}A_{s,j}^{}A_{i,j}^{}A_{s,j}^{-1}A_{r,j}^{-1}&\text{if $r=i<s<j$,}\\
[A_{r,j}^{},A_{s,j}^{}]A_{i,j}^{}[A_{r,j}^{},A_{s,j}^{}]^{-1}&\text{if $r<i<s<j$,}
\end{cases}
\end{equation}
where $[u,v]=uvu^{-1}v^{-1}$ denotes the commutator.  
See, for instance, Birman \cite{Bir75} as a general reference on braid groups. 
It is well known that the pure braid group admits a direct product decomposition 
\begin{equation} \label{eq:directprod}
P_n=Z \times P_n/Z,
\end{equation} 
where $Z=Z(P_n) \cong \Z$ is the center of $P_n$, generated by 
\begin{equation} \label{eq:center}
Z_n=(A_{1,2})(A_{1,3}A_{2,3})\cdots (A_{1,n}\cdots A_{n-1,n}).
\end{equation} 
Note that $P_3=Z(P_3) \times P_3/Z(P_3) \cong \Z \times F_2$, where $F_2$ is the free group on two generators. 
For any $n\ge 3$, by (\ref{eq:directprod}),
there is a split, short exact sequence
\begin{equation} \label{eq:autseq}
1\to \tv(P_n) \to \Aut(P_n) \leftrightarrows \Aut(P_n/Z) \to 1,
\end{equation}
where the subgroup $\tv(P_n)$ of $\Aut(P_n)$ consists of those automorphisms which become trivial upon passing to the quotient $P_n/Z$.  

For a group $G$ with infinite cyclic center $Z=Z(G)=\langle z\rangle$, a transvection is an endomorphism of $G$ of the form $x \mapsto xz^{t(x)}$, where $t\colon G \to \Z$ is a homomorphism, see Charney and Crisp \cite{CC05}.  Such a map is an automorphism if and only if its restriction to $Z$ is surjective, which is the case if and only if $z \mapsto z$ or $z \mapsto z^{-1}$, that is, $t(z)=0$ or $t(z)=-2$.  For the pure braid group $P_n$, the transvection subgroup $\tv(P_n)$ of $\Aut(P_n)$ consists of automorphisms of the form $A_{i,j}\mapsto A_{i,j}Z_n^{t_{i,j}}$, where $t_{i,j}\in\Z$ and $\sum t_{i,j}$ is either equal to $0$ or $-2$. In the former case, $Z_n \mapsto Z_n$, while $Z_n \mapsto Z_n^{-1}$ in the latter.  This yields a surjection $\tv(P_n) \to \Z_2$, with kernel consisting of transvections for which $\sum t_{i,j}=0$.  Since $P_n$ has $\binom{n}{2}=N+1$ generators, this kernel is free abelian of rank $N$.  The choice $t_{1,2}=-2$ and all other $t_{i,j}=0$ gives a splitting $\Z_2 \to \tv(P_n)$.  Thus, 
$\tv(P_n) \cong \Z^N \rtimes \Z_2$.  Explicit generators of $\tv(P_n)$ are given below.

For $n\ge 4$, Bell and Margalit \cite{BM07} show that the automorphism group of the pure braid group admits a semidirect product decomposition
\begin{equation} \label{eq:semi}
\Aut(P_n) \cong (\Z^N \rtimes \Z_2) \rtimes \Mod(\bS_{n+1}).
\end{equation}
Here, $\Z^N \rtimes \Z_2=\tv(P_n)$ is the transvection subgroup of $\Aut(P_n)$ described above, $\bS_{n+1}$ denotes the sphere $S^2$ with $n+1$ punctures, and 
$\Mod(\bS_{n+1})$ is the extended mapping class group of $\bS_{n+1}$, the group of isotopy classes of all self-diffeomorphisms of $\bS_{n+1}$.
The semidirect product decomposition \eqref{eq:semi} is used in \cite{Coh10} to determine a finite presentation for $\Aut(P_n)$.  From this work, it follows that $\Aut(P_n)$ is generated by automorphisms
\begin{equation} \label{eq:autgens}
\xi, \ \beta_k\ (1\le k\le n),\ \psi,\ \phi_{p,q}\  (1\le p < q \le n,\ \{p,q\}\neq \{1,2\}),
\end{equation}
given explicitly by
\begin{equation} \label{eq:autos}
\begin{aligned}
\xi\colon A_{i,j}&\mapsto (A_{i+1,j} \cdots A_{j-1,j})^{-1} A_{i,j}^{-1} (A_{i+1,j} \cdots A_{j-1,j}),\\
\beta_k\colon A_{i,j}&\mapsto \begin{cases}
A_{i-1,j}&\text{if $k=i-1$,}\\
A_{i,i+1}^{-1}A_{i+1,j}^{}A_{i,i+1}^{}&\text{if $k=i<j-1$,}\\
A_{i,j-1}&\text{if $k=j-1>i$,}\\
A_{j,j+1}^{-1}A_{i,j+1}^{}A_{j,j+1}^{}&\text{if $k=j$,}\\
A_{i,j}&\text{otherwise,}
\end{cases}
\quad \text{for $1\le k \le n-1$,}\\
\beta_n\colon A_{i,j}&\mapsto \begin{cases}
A_{i,j}&\text{if $j \neq n$,}\\
(A_{i,n}A_{1,i}\cdots A_{i-1,i} A_{i,i+1}\cdots A_{i,n-1})^{-1}&\text{if  $j=n$,}
\end{cases}\\
\psi\colon A_{i,j}& \mapsto \begin{cases} A_{1,2}Z_n^{-2}&\text{if $i=1$ and $j=2$,}\\ A_{i,j}&\text{otherwise,} \end{cases} \\
\phi_{p,q}\colon A_{i,j}& \mapsto \begin{cases} A_{1,2}Z_n^{}&\text{if $i=1$ and $j=2$,}\\ A_{p,q}Z_n^{-1}&\text{if $i=p$ and $j=q$,}\\ A_{i,j}&\text{otherwise.} \end{cases}
\end{aligned}
\end{equation}
It is readily checked that these are all automorphisms of $P_n$.  The automorphisms $\psi$ and $\phi_{p,q}$ are transvections. 
For
$k\le n-1$, 
$\beta_k\in\Aut(P_n)$ 
arises from the conjugation action of $B_n$ on $P_n$, $\beta_k(A_{i,j})=\sigma_k^{-1}A_{i,j}^{}\sigma_k^{}$, see Dyer and Grossman \cite{DG81}.

\begin{remark}
The presentation of $\Aut(P_n)$ found in \cite{Coh10} is given in terms of the generating set $\epsilon$, $\omega_k$ ($1\le k \le n$), $\psi$, $\phi_{p,q}$ ($1\le p < q \le n,\ \{p,q\}\neq \{1,2\}$), where $\xi=\epsilon \circ\psi$, $\beta_2=\omega_2\circ \phi_{1,3}^{-1}$, $\beta_n=\omega_n\circ\psi\circ\phi_{1,n}\circ\phi_{2,n}$, and $\beta_k=\omega_k$ for $k\neq 2,n$.  This presentation exhibits the semidirect product structure \eqref{eq:semi} of $\Aut(P_n)$.
\end{remark}

\section{Epimorphisms to free groups}

We study surjective homomorphisms from the pure braid group $P_n$ to the free group $F_k$ on $k\ge 2$ generators.  Since $P_2=\Z$ is infinite cyclic, we assume that $n\ge 3$.  We begin by exhibiting a number of specific such homomorphisms.

Let $F_2=\langle x,y\rangle$ be the free group on two generators, and write $[n]=\{1,2,\dots,n\}$.  
If $I=\{i,j,k\} \subset [n]$ with $i<j<k$, define $f_I\colon P_n \to F_2$ by
\begin{equation} \label{eq:three}
f_I(A_{r,s})=\begin{cases}
x&\text{if $r=i$ and $s=j$,}\\
y&\text{if $r=i$ and $s=k$,}\\
y^{-1}x^{-1}&\text{if $r=j$ and $s=k$,}\\
1&\text{otherwise.}
\end{cases}
\end{equation}
If $I=\{i,j,k,l\} \subset [n]$ with $i<j<k<l$, define $f_I\colon P_n \to F_2$ by
\begin{equation} \label{eq:four}
f_I(A_{r,s})=\begin{cases}
x&\text{if $r=i$ and $s=j$,}\\
y&\text{if $r=i$ and $s=k$,}\\
y^{-1}x^{-1}&\text{if $r=j$ and $s=k$,}\\
y^{-1}x^{-1}&\text{if $r=i$ and $s=l$,}\\
x^{}y^{}x^{-1}&\text{if $r=j$ and $s=l$,}\\
x&\text{if $r=k$ and $s=l$,}\\
1&\text{otherwise.}
\end{cases}
\end{equation}
In either case \eqref{eq:three} or \eqref{eq:four}, note that $f_I$ is surjective by construction.  
It is readily checked that $f_I$ is a homomorphism.
We will show that these are, in an appropriate sense, the only epimorphisms from the pure braid group to a nonabelian free group.

\begin{remark} \label{rem:topological}
The epimorphisms $f_I \colon P_n \to F_2$ are induced by maps of topological spaces.  Let 
\[F(\C,n)=\{(z_1,\dots,z_n) \in \C^n \mid z_i \neq z_j\ \text{if}\ i \neq j\}\]
be the configuration space of $n$ distinct ordered points in $\C$.  It is well known that $P_n=\pi_1(F(\C,n))$ and that $F(\C,n)$ is a $K(P_n,1)$-space.

For a subset $I$ of $[n]$ of cardinality $k$, let $p_I \colon F(\C,n) \to F(\C,k)$ denote the projection which forgets all coordinates not indexed by $I$.  The induced map on pure braid groups forgets the corresponding strands.  Additionally, let $q_n \colon F(\C,n) \to F(\C,n)/\C^*$ denote the natural projection, where $\C^*$ acts by scalar multiplication.  
In particular, $q_3 \colon F(\C,3) \to F(\C,3)/\C^*\cong \C \times (\C\smallsetminus\{\text{two points}\})$.  
Finally, define $g \colon F(\C,4) \to F(\C,3)$ by
\[
g(z_1,z_2,z_3,z_4) = \bigl((z_1+z_2-z_3-z_4)^2,(z_1+z_3-z_2-z_4)^2,(z_1+z_4-z_2-z_3)^2\bigr).
\]
One can check that if $|I|=3$, then $f_I = (q_3 \circ p_I)_*$, while if $|I|=4$, $f_I=(q_3 \circ g \circ p_I)_*$. 

In the case $n=4,$ $F(\C,n)$ is diffeomorphic to $\C \times M$, where $M$ is the complement of the Coxeter arrangement \A\ of type $D_3$. With this identification, the mappings $p_I$ and $g$ correspond to the components of the mapping constructed in \cite[Example 3.25]{CFR10}. The map $g$ is the pencil associated with the non-local component of the first resonance variety of $H^*(M;\C)$ as in \cite{LY00,FY}, see below. 

One can also check that the homomorphism $g_*\colon P_4 \to P_3$ is the restriction to pure braid groups of the famous homomorphism $B_4 \to B_3$ of full braid groups, given by $\sigma_1 \mapsto \sigma_1$, $\sigma_2 \mapsto \sigma_2$, $\sigma_3 \mapsto \sigma_1$.
\end{remark}

Let $G$ and $H$ be groups, and let $f$ and $g$ be (surjective) homomorphisms from $G$ to $H$.  Call $f$ and $g$ equivalent if there are automorphisms $\phi\in \Aut(G)$ and $\psi \in \Aut(H)$ so that $g \circ \phi = \psi \circ f$.  If $f$ and $g$ are equivalent, we write $f\sim g$. 

\begin{proposition} \label{prop:34equiv}
If $I$ and $J$ are subsets of $[n]$ of cardinalities $3$ or $4$, then the epimorphisms $f_I$ and $f_J$ from $P_n$ to $F_2$ are equivalent.
\end{proposition}
\begin{proof}
If $n=3$, there is nothing to prove. So assume that $n\ge 4$.

Let $I=\{i_1,\dots,i_q\} \subset [n]$ with $q \ge 2$ and $i_1<i_2<\dots<i_q$.  Define $\alpha_I \in B_n$ by
\[
\alpha_I = (\sigma_{i_1-1} \cdots \sigma_1)(\sigma_{i_2-1} \cdots \sigma_2) \cdots (\sigma_{i_q-1} \cdots \sigma_q),
\]
where $\sigma_{i_k-1} \cdots \sigma_k=1$ if $i_k=k$. 
Then $\alpha_I^{-1}A_{i_r,i_s}^{}\alpha_I^{}=A_{r,s}^{}$ for $1\le r < s \le q$.  This can be seen by checking that, for instance, the geometric braids $\alpha_I^{}A_{r,s}^{}\alpha_I^{-1}$ and $A_{i_r,i_s}$ are equivalent.  Denote the automorphism $A_{i,j}\mapsto 
\alpha_I^{-1}A_{i,j}^{}\alpha_I^{}$ of $P_n$ by the same symbol,
\[
\alpha_I=(\beta_{i_1-1} \cdots \beta_1)(\beta_{i_2-1} \cdots \beta_2) \cdots (\beta_{i_q-1} \cdots \beta_q)
\in\Aut(P_n).
\]
Then, $\alpha_I(A_{i_r,i_s})=A_{r,s}$ for 
$1\le r < s \le q$.

If $I$ has cardinality $3$, then, by the above, we have $f_I = f_{[3]} \circ \alpha_I$, so $f_I \sim f_{[3]}$.  Similarly, if $|I|=4$, then $f_I = f_{[4]} \circ \alpha_I$ and $f_I \sim f_{[4]}$.  Thus it suffices to show that $f_{I} \sim f_{[3]}$ for some $I$ with $|I|=4$.  This can be established by checking that 
$f_{I}=f_{[3]} \circ \beta_n$ for $I=\{1,2,3,n\}$.
\end{proof}

\begin{remark} The homomorphisms $f_I$ also have a natural interpretation in terms of the moduli space $\MM_{0,n}$ of genus-zero curves with $n$ marked points. By definition, $\MM_{0,n}$ is the quotient of the configuration space $F(S^2,n)$ of the Riemann sphere $S^2=\C\cup \{\infty\}$ by the action of $\PSL(2,\C)$. The map $h_n \colon F(\C,n) \to \MM_{0,n+1}$ given by $h_n(z_1, \ldots, z_n) = [(z_1, \ldots, z_n, \infty)]$ induces a homeomorphism 
\[
F(\C,n)/\Aff(\C) \to \MM_{0,n+1},
\] where $\Aff(\C)\cong \C \rtimes \C^*$ is the affine group, and a homotopy equivalence 
\[
\bar{h}_n \colon F(\C,n)/\C^* \to \MM_{0,n+1}.
\] 

For $1\leq i \leq 5$, let $\delta_i \colon \MM_{0,5} \to \MM_{0,4}$ be defined by forgetting the $i$-th point. Then, in the notation of Remark~\ref{rem:topological}, for $1\leq i \leq 4$, $\delta_i\circ h_4= \bar{h}_3\circ q_3\circ p_I$, where $I=[4]\setminus \{i\}$. Up to a linear change of coordinates in $F(\C,3)$, $\delta_5 \circ h_4=\bar{h}_3 \circ q_3 \circ g$. (See also Pereira \cite[Example 3.1]{Per10}.) The maps $\delta_i \colon \MM_{0,5} \to \MM_{0,4}$, $1\leq i \leq 5$, are clearly equivalent up to diffeomorphism of the source. Applying Remark~\ref{rem:topological}, this gives an alternate proof of Proposition~\ref{prop:34equiv}.
\end{remark}

For our next result, we require some properties of the cohomology ring of the pure braid group. 
Let $A=\bigoplus_{k\geq 0} A^k$ be a 
connected, graded-commutative 
algebra over a 
field $\k$, with $\dim A^1< \infty$. 
Since $a\cdot a=0$ for each $a\in A^1$, multiplication 
by $a$ defines a cochain complex $(A,\delta_a)$:
\[
A^0\xrightarrow{\quad \delta_a\quad} A^1 \xrightarrow{\quad \delta_a\quad} A^2 \xrightarrow{\quad \delta_a\quad} \cdots\cdots \xrightarrow{\quad \delta_a\quad} A^\ell,
\]
where $\delta_a(x)=ax$.  The resonance varieties $\RR^d(A)$ of $A$ are defined by
\[
\RR^d(A) = \{a \in A^1 \mid H^d(A,\delta_a) \neq 0\}.
\]
If $\dim A^d<\infty$, then $\RR^d(A)$ is an algebraic set in $A^1$.

In the case where $A=H^*(M(\A);\k)$ is the cohomology ring of the complement of a complex hyperplane arrangement, and \k\ has characteristic zero, work of Libgober and Yuzvinsky \cite{LY00} (see also \cite{FY}) shows that $\RR^1(A)$ is the union of the maximal isotropic subspaces of $A^1$ for the quadratic form
\begin{equation} \label{eq:mult}
\mu\colon A^1 \otimes A^1 \to A^2, \ \mu(a\otimes b) = a b
\end{equation}
having dimension at least two. 
Note that, for any field \k, any isotropic subspace of $A^1$ of dimension at least two is contained in $\RR^1(A)$.

For our purposes, it will suffice to take $\k=\Q$.  Let $A=H^*(P_n;\Q)$ be the rational cohomology ring of the pure braid group, that is, the cohomology of the complement of the braid arrangement in $\C^n$.  By work of Arnold \cite{Arn} and F.~Cohen \cite{FC}, $A$ is generated by degree one elements $a_{i,j}$, $1\le i<j\le  n$, which satisfy (only) the relations
\[
a_{i,j}a_{i,k}-a_{i,j}a_{j,k}+a_{i,k}a_{j,k}\ \text{for}\ 1\le i<j<k\le n,
\]
and their consequences.  The irreducible components of the first resonance variety $\RR^1(A)$ may be obtained from work of Cohen and Suciu \cite{CS} 
(see also \cite{Per10}), 
and may be described as follows.  The rational vector space $A^1=H^1(P_n;\Q)$ is of dimension $\binom{n}{2}$, and has basis $\{a_{i,j}
\mid 1\le i < j\le n\}$.  If $I=\{i,j,k\}\subset [n]$ with $i<j<k$, let $V_I$ be the subspace of $A^1$ defined by
\[
V_I=\Span\{a_{i,j}-a_{j,k}, a_{i,k}-a_{j,k}\}.
\]
If $I=\{i,j,k,l\}\subset [n]$ with $i<j<k<l$, let $V_I$ be the subspace of $A^1$ defined by
\[
V_I=\Span\{a_{i,j}+a_{k,l}-a_{j,k}-a_{i,l}, a_{i,k}+a_{j,l}-a_{j,k}-a_{i,l}\}.
\]
The $2$-dimensional subspaces $V_I$ of $A^1$, where $|I|=3$ or $|I|=4$, are the irreducible components of $\RR^1(A)=\bigcup_{k=3}^4\bigcup_{|I|=k}V_I$.

\begin{remark}
One can check that $V_I=f_I^*(H^1(F_2;\Q))$, where $f_I\colon P_n \to F_2$ is defined by \eqref{eq:three} if $|I|=3$ and by \eqref{eq:four} if $|I|=4$. Work of Schenck and Suciu \cite[Lemma 5.3]{SchSuc2} implies that for any two components $V_I$ and $V_J$ of $\RR^1(A)$, there is an isomorphism of $A^1=H^1(P_n;\Q)$ taking $V_I$ to $V_J$.  This provides an analog, on the level of (degree-one) cohomology, of Proposition \ref{prop:34equiv}.
\end{remark}

\begin{theorem} \label{thm:PntoFk}
Let $n\ge 3$ and $k \ge 2$, and consider the pure braid group $P_n$ and the free group $F_k$.
\begin{enumerate}
\item \label{item:one}
If $k \ge 3$, there are no epimorphisms from $P_n$ to $F_k$.
\item \label{item:two}
If $k=2$, there is a single equivalence class of epimorphisms from $P_n$ to~$F_2$.
\end{enumerate}
\end{theorem}
\begin{proof}
For part \eqref{item:one}, if $f\colon P_n \to F_k$ is an epimorphism, then 
$f$ splits, so 
\[
f^*\colon H^1(F_k;\Q) 
\to H^1(P_n;\Q)
\]
 is injective. Consequently, 
$f^*(H^1(F_k;\Q))$ is a $k$-dimensional isotropic subspace of $A^1=H^1(P_n;\Q)$ for the form \eqref{eq:mult}.  Since this subspace is isotropic, it must be contained in an irreducible component of $\RR^1(A)$.  Since these components are all of dimension $2$, we cannot have $k\ge 3$.

For part \eqref{item:two}, by Proposition \ref{prop:34equiv}, it suffices to show that an arbitrary epimorphism $f\colon P_n\to F_2$ is equivalent to $f_I$ for some $I$ of cardinality $3$ or $4$.  We will extensively use that fact that if $[a,b]=1$ in $F_2$, then $\langle a,b\rangle<F_2$ is free and abelian, so $a=z^m$ and $b=z^n$ for some $z\in F_2$ and $m,n\in\Z$.  
Additionally, if the homology class $[a]$ of $a$ in $H_1(F_2;\Z)$ is part of a basis and $[a,b]=1$, then $b$ is a power of $a$. Indeed, suppose  $\{[a],[a']\}$ is a basis for $H_1(F_2;\Z)$. Write $[z]=c_1[a]+c_2[a']$ with $c_1,c_2 \in \Z$.  Then we have $[a]=m[z]=mc_1[a]+mc_2[a']$, which implies $c_2=0$ and $m=c_1=\pm 1$, yielding the assertion.

Let $f\colon P_n\to F_2$ be an epimorphism.  
Then $f^*(H^1(F_2;\Q))$ is a $2$-dimensional isotropic subspace of $H^1(P_n;\Q)$ for the form \eqref{eq:mult}.  Since the irreducible components of $\RR^1(A)$ are all of dimension $2$, we must have $f^*(H^1(F_2;\Q))=V_I$ for some $I\subset[n]$ of cardinality $3$ or $4$.  Since the cohomology rings of  $P_n$ and $F_2$ are torsion-free, passing to integer coefficients, we have $f^*(H^1(F_2;\Z))=
V_I \cap \Z^{\binom{n}{2}}
\subset H^1(P_n;\Z)$.  Consequently, there is an automorphism $\varphi^*$ of $H^1(F_2;\Z)$ so that $f^*=f_I^*\circ \varphi^*$.  Passing to homology (again using torsion-freeness), we have $f_*=\varphi_* \circ (f_I)_*$, where $\varphi_*\in\Aut(H_1(F_2;\Z))$ is dual to $\varphi^*$.  Let $\varphi\in\Aut(F_2)$ be an automorphism which induces $\varphi_*$.

 From the definitions \eqref{eq:three} and \eqref{eq:four} of the epimorphisms $f_I$, there exists
$\{i,j,k\}$ with $1\le i<j<k\le n$, $f_I(A_{i,j})=x$, and $f_I(A_{i,k})=y$, where $F_2=\langle x,y\rangle$.  Let $u=f(A_{i,j})$ and $v=f(A_{i,k})$.  
Using the equation $f_*=\varphi_* \circ (f_I)_*$, we have 
\[
[u]= [f(A_{i,j})] = f_*([A_{i,j}]) = \varphi_*([f_I(A_{i,j})] = \varphi_*([x]) = [\varphi(x)],
\]
and similarly $[v]=[\varphi(y)]$.  Thus $\{[u],[v]\}$ is a basis for $H_1(F_2;\Z)$.

Let $w=f(A_{j,k})$.  Using the pure braid relations \eqref{eq:purebraidrels}, we have $A_{i,j}A_{i,k}A_{j,k}=A_{i,k}A_{j,k}A_{i,j}=A_{j,k}A_{i,j}A_{i,k}$.  Applying $f$, these imply that $[uv,w]=1$ and $[u,vw]=1$ in $F_2$. 
Since $\{[u],[u]+[v]\}$ is a basis for $H_1(F_2;\Z)$, these imply that $w=(uv)^m$ and $vw=u^n$ for some $m,n\in \Z$.  A calculation with homology classes reveals that $m=n=-1$.  Hence $w=v^{-1}u^{-1}$, i.e., $uvw=1$.

Suppose that $f(A_{r,s})=1$ for all $\{r,s\}\not\subset\{i,j,k\}$.  Then the image of $f$ is contained in the subgroup $\langle u,v\rangle$ of $F_2$.  Since $f$ is by hypothesis an epimorphism, we have $\langle u,v\rangle=F_2$.  Letting $\lambda$ be the automorphism of $F_2$ taking $u$ to $x$ and $v$ to $y$, we have $\lambda \circ f = f_{\{i,j,k\}}$.

Now suppose that $f(A_{r,s})\neq 1$ for some $\{r,s\}\not\subset\{i,j,k\}$. First assume that $\{r,s\} \cap \{i,j,k\}=\emptyset$.  We claim that $f(A_{r,s})=1$. There are various cases depending on the relative positions of $r<s$ and $i<j<k$.  We consider the case $i<r<j<k<s$ and leave the remaining analogous cases to the reader.  In this instance, we have relations $[A_{j,k},A_{r,s}]=1$ and $A_{i,j}^{-1}A_{r,s}^{}A_{i,j}^{}=[A_{i,s}^{},A_{j,s}^{}]A_{r,s}^{}[A_{i,s}^{},A_{j,s}^{}]^{-1}$.  The second of these, together with the pure braid relations \eqref{eq:purebraidrels} may be used to show that $[A_{i,j},A_{j,s}^{-1}A_{r,s}^{}A_{j,s}^{}]=1$.  Applying $f$, we have $[w,f(A_{r,s})]=1$ and $[u,z^{-1}f(A_{r,s})z]=1$ in $F_2$, where $z=f(A_{j,s})$.  Since any two element subset of $\{[u],[v],[w]\}$ forms a basis for $H_1(F_2;\Z)$, these relations imply that $f(A_{r,s})=w^m$ and $z^{-1}f(A_{r,s})z=u^n$ for some $m,n\in \Z$.  Consequently, $m[w]=n[u]$ in $H_1(F_2;\Z)$, which forces $m=n=0$ and $f(A_{r,s})=1$.

Thus, we must have $f(A_{r,s})\neq 1$ for some $\{r,s\}$ with $|\{r,s\} \cap \{i,j,k\}|=1$.  As above, there are several cases, and we consider a representative one, leaving the other, similar, cases to the reader.

Assume that $r=k$, so that $i<j<k<s$, and that $f(A_{k,s})\neq 1$.  Applying $f$ to the pure braid relations $[A_{i,j},A_{k,s}]=1$, $[A_{j,k},A_{i,s}]=1$, and $[A_{i,k}^{},A_{k,s}^{-1}A_{j,s}^{}A_{k,s}^{}]=1$ yields $[u,f(A_{k,s})]=1$, $[w,f(A_{i,s})]=1$, and $[v,f(A_{k,s}^{-1}A_{j,s}^{}A_{k,s}^{})]=1$ in $F_2$. 
It follows that $f(A_{k,s})=u^m$, $f(A_{i,s})=w^n$, and $f(A_{j,s})=u^mv^lu^{-m}$ for some $m,n,l\in \Z$.  Since $f(A_{k,s})\neq 1$, we have $m\neq 0$.  Then, applying $f$ to the pure braid relations
$[A_{i,k},A_{i,s}A_{k,s}] = 1$ and $[A_{j,k},A_{j,s}A_{k,s}] = 1$,
we obtain $[v,w^nu^m]=1$ and $[w,u^mv^l]=1$ in $F_2$.  Thus, $w^nu^m=v^p$ and $u^mv^l=w^q$ for some $p,q\in\Z$.  Passing to homology, using the fact that $[u]+[v]+[w]=0$ in $H_1(F_2;\Z)$ since $uvw=1$ in $F_2$, reveals that $m=n=l$.

In $P_n$, we also have the relation $[A_{j,s},A_{k,s}A_{j,k}]=1$.  Applying $f$ we obtain the relation $[u^mv^mu^{-m},u^mw]=1$ in $F_2$.  Hence, $u^mv^mu^{-m}=z^p$ and $u^mw=z^q$ for some $z\in F_2$ and $p,q\in\Z$.  Writing $[z]=c_1[u]+c_2[v]$, we have
\[
m[v] = pc_1[u]+pc_2[v] \ \text{and}\ (m-1)[u] - [v] = qc_1[u]+qc_2[v]
\]
in $H_1(F_2;\Z)$.  It follows that $m=n=l=1$, and therefore $f(A_{i,s})=w=v^{-1}u^{-1}$, $f(A_{j,s})=uvu^{-1}$, and $f(A_{k,s})=u$.
Thus the image of $f$ is contained in the subgroup $\langle u,v\rangle$ of $F_2$.  As before, $\langle u,v\rangle=F_2$ since $f$ is an epimorphism.  Recalling that $\lambda$ is the automorphism of $F_2$ taking $u$ to $x$ and $v$ to $y$, we have $\lambda \circ f = f_{\{i,j,k,s\}}$.
\end{proof}

The proof of Theorem~\ref{thm:PntoFk}(a) actually yields the following more general result. 
\begin{theorem} \label{thm:resonance}
Let $G$ be a finitely generated group, and $\k$ an algebraically closed field. Then there are no epimorphisms from $G$ to $F_k$ for $k > \dim \RR^1(H^*(G,\k))$.
\end{theorem}

Recall that the corank of a group $G$ is the largest natural number $k$ for which the free group $F_k$ is an epimorphic image of $G$.  The corank of $P_2\cong\Z$ is $1$. For larger $n$, as an immediate consequence of Theorem \ref{thm:PntoFk}, we obtain the following. 

\begin{corollary} \label{cor:corank}
For $n \ge 3$, the corank of the pure braid group $P_n$ is equal to $2$.
\end{corollary}

Theorem \ref{thm:resonance} yields a more general result.

\begin{corollary}\label{cor:upper}
 Let $G$ be a finitely generated group and \k\ an algebraically closed field. Then the corank of $G$ is bounded above by $\dim \RR^1(H^*(G,\k))$.
\end{corollary}

In fact, results of Dimca, Papadima, and Suciu \cite{DPS09} imply that the corank of $G$ is equal to $\dim \RR^1(H^*(G,\C))$ for a wide class of quasi-Kahler groups, including fundamental groups of complex projective hypersurface complements.

\section{Pure braid groups are not residually free}

Recall that a group $G$ is residually free if for every $x \neq 1$ in $G$, there is a homomorphism $f$ from $G$ to a free group $F$ so that $f(x)\neq 1$ in $F$.  In this section, we show that $P_n$ is not residually free for $n\ge 4$, and derive some consequences.  Since $P_2\cong\Z$ and $P_3\cong \Z\times F_2$, these groups are residually free.  

For each $I\subset [n]$ of cardinality $3$ or $4$, let $K_I=\ker(f_I\colon P_n \to F_2)$.  Define
\[
K_n = \bigcap_{k=3}^4\bigcap_{|I|=k} K_I.
\]

\begin{proposition} \label{prop:char}
The subgroup $K_n$ of $P_n$ is characteristic.
\end{proposition}
\begin{proof}
It suffices to show that if $\beta$ is one of the generators of $\Aut(P_n)$ listed in \eqref{eq:autgens}, then $\beta(K_n)=K_n$.  

If $\beta=\psi$ or $\beta=\phi_{p,q}$ is a transvection, then for each $I$, $f_I \circ \beta= f_I$ since $f_I(Z_n)=1$, where $Z_n$ is the generator of the center $Z(P_n)$ of $P_n$ recorded in \eqref{eq:center}.  It follows that $\beta(K_I)=K_I$ for each $I$, which implies that $\beta(K_n)=K_n$.

If $\beta=\xi$, then for each $I$, it is readily checked that $f_I \circ \xi = \lambda \circ f_I$, where $\lambda \in \Aut(F_2)$ is defined by $\lambda(x)=x^{-1}$ and $\lambda(y)=xy^{-1}x^{-1}$.  Thus, $\xi(K_I)=K_I$, and $\xi(K_n)=K_n$.

If $\beta=\beta_k$ for $1\le k \le n-1$, let $\tau_k$ denote the permutation induced by $\beta_k$.  Let $I=\{i_1,\dots,i_l\}$ where $l=3$ or $4$, and let $\tau_k(I)$ denote the set $\{\tau_k(i_1),\dots,\tau_k(i_l)\}$ with the elements in increasing order.  Define automorphisms $\lambda_1,\lambda_2\in\Aut(F_2)$ by
\[
\lambda_1\colon\begin{cases}x \mapsto x,\\y\mapsto x^{-1}y^{-1},\end{cases}\quad\text{and}\quad
\lambda_2\colon\begin{cases}x \mapsto xyx^{-1},\\y\mapsto x,\end{cases}
\]
and set $\lambda_3=\lambda_1$.  Then, calculations with the definitions of the automorphism $\beta_k$ and the epimorphisms $f_I\colon P_n \to F_2$ (see \eqref{eq:autos}, \eqref{eq:three}, and \eqref{eq:four}) reveal that
\[
f_I \circ \beta_k = \begin{cases}
\lambda_j \circ f_{\tau_k(I)}&\text{if $k=i_j = i_{j+1}-1$,}\\
f_{\tau_k(I)}&\text{otherwise.}
\end{cases}
\]
Note that, in the first case above, $\tau_k(I)=I$ and $j<l$.  Thus, $K_I=\beta_k(K_{\tau_k(I)})$ for each $I$, and $\beta_k$ permutes the subgroups $K_I$ of $P_n$ (for $|I|=3$ and $|I|=4$ respectively).  It follows that $\beta_k(K_n)=K_n$.

Finally, if $\beta=\beta_n$, calculations with \eqref{eq:autos}, \eqref{eq:three}, and \eqref{eq:four}) reveal that
\[
f_I \circ \beta_n = \begin{cases}
f_{I\cup\{n\}}&\text{if $|I|=3$ and $n\notin I$,}\\
\lambda_1 \circ f_{I}&\text{if $|I|=3$ and $n\in I$,}\\
f_{I}&\text{if $|I|=4$ and $n\notin I$,}\\
f_{I\smallsetminus\{n\}}&\text{if $|I|=4$ and $n\in I$.}
\end{cases}
\]
Thus, $\beta_n(K_I)=K_I$ if either $|I|=3$ and $n\in I$ or $|I|=4$ and $n\notin I$, while $\beta_n(K_I)=K_{I\cup\{n\}}$ and $\beta_n(K_{I\cup\{n\}})=K_I$ if $|I|=3$ and $n\notin I$.  It follows that $\beta_n(K_n)=K_n$, and $K_n\ \text{char}\ P_n$.
\end{proof}

\begin{theorem} \label{thm:notresfree}
For $n\ge 4$, the pure braid group $P_n$ is not residually free.
\end{theorem}
\begin{proof}
Let $f \colon P_n \to F_k$ be a surjective homomorphism from the pure braid group to a nonabelian free group.  We claim that $K_n$ is contained in the kernel of $f$.  By Theorem \ref{thm:PntoFk},  $k= 2$ and $f \sim f_{[3]}$.  Thus there are automorphisms $\alpha \in \Aut(P_n)$ and $\lambda \in \Aut(F_2)$ so that $\lambda \circ f = f_{[3]} \circ \alpha$.  Let $x \in K_n$.  Then $\alpha(x) \in K_n$ since $K_n$ is characteristic in $P_n$ by Proposition \ref{prop:char}.  Since $K_n \subset K_{[3]}=\ker(f_{[3]})$ by definition, we have $\alpha(x) \in \ker(f_{[3]})$. Hence, $\lambda \circ f(x)=f_{[3]} \circ\alpha(x)=1$, and $x \in \ker(f)$.

To complete the proof, it suffices to exhibit a nontrivial element of $K_n$ that is in the kernel of every homomorphism $g\colon P_n \to \Z$ from the pure braid group to an abelian free group.  This is straightforward since is it easy to see that the intersection $K_n \cap [P_n,P_n]$ of $K_n$ with the commutator subgroup of $P_n$ is nontrivial.  For instance, a calculation reveals that $x=[[A_{1,2},A_{2,3}],[A_{2,3},A_{3,4}]] \in K_n \cap [P_n,P_n]$.  
The pure braid $x$ is nontrivial (one can check that the braids $[A_{1,2},A_{2,3}][A_{2,3},A_{3,4}]$ and $[A_{2,3},A_{3,4}][A_{1,2},A_{2,3}]$ are distinguished by the Artin representation).  
Thus $x\neq 1$ is in the kernel of every homomorphism from $P_n$ to a free group, and $P_n$ is not residually free.
\end{proof}

\begin{figure}
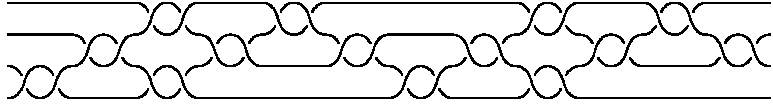
  
\begin{center}
\ifpic \leavevmode 
\xygraph{ !{0;/r1pc/:}
!{\xcaph[1]@(0)}[ld] 
!{\xcaph[1]@(0)}[ld] 
!{\htwistneg}[u(2)] 
!{\xcaph[1]@(0)}[ld] 
!{\xcaph[1]@(0)}[ld] 
!{\htwistneg}[u(2)] 
!{\xcaph[1]@(0)}[ld] 
!{\htwistneg}[l(1)d(2)] 
!{\xcaph[1]@(0)}[u(3)] 
!{\xcaph[1]@(0)}[ld] 
!{\htwistneg}[l(1)d(2)] 
!{\xcaph[1]@(0)}[u(3)] 
!{\htwistneg}[l(1)d(2)] 
!{\htwist}[u(2)] 
!{\htwistneg}[l(1)d(2)] 
!{\htwist}[u(2)] 
!{\xcaph[1]@(0)}[ld] 
!{\htwist}[l(1)d(2)] 
!{\xcaph[1]@(0)}[u(3)] 
!{\xcaph[1]@(0)}[ld] 
!{\htwist}[l(1)d(2)] 
!{\xcaph[1]@(0)}[u(3)] 
!{\htwist}[l(1)d(2)] 
!{\xcaph[1]@(0)}[ld] 
!{\xcaph[1]@(0)}[u(3)] 
!{\htwist}[l(1)d(2)] 
!{\xcaph[1]@(0)}[ld] 
!{\xcaph[1]@(0)}[u(3)] 
!{\xcaph[1]@(0)}[ld] 
!{\htwistneg}[l(1)d(2)] 
!{\xcaph[1]@(0)}[u(3)] 
!{\xcaph[1]@(0)}[ld] 
!{\htwistneg}[l(1)d(2)] 
!{\xcaph[1]@(0)}[u(3)] 
!{\xcaph[1]@(0)}[ld] 
!{\xcaph[1]@(0)}[ld] 
!{\htwistneg}[u(2)] 
!{\xcaph[1]@(0)}[ld] 
!{\xcaph[1]@(0)}[ld] 
!{\htwistneg}[u(2)] 
!{\xcaph[1]@(0)}[ld] 
!{\htwist}[l(1)d(2)] 
!{\xcaph[1]@(0)}[u(3)] 
!{\xcaph[1]@(0)}[ld] 
!{\htwist}[l(1)d(2)] 
!{\xcaph[1]@(0)}[u(3)] 
!{\htwistneg}[l(1)d(2)] 
!{\htwist}[u(2)] 
!{\htwistneg}[l(1)d(2)] 
!{\htwist}[u(2)] 
!{\xcaph[1]@(0)}[ld] 
!{\htwistneg}[l(1)d(2)] 
!{\xcaph[1]@(0)}[u(3)] 
!{\xcaph[1]@(0)}[ld] 
!{\htwistneg}[l(1)d(2)] 
!{\xcaph[1]@(0)}[u(3)] 
!{\htwist}[l(1)d(2)] 
!{\xcaph[1]@(0)}[ld] 
!{\xcaph[1]@(0)}[u(3)] 
!{\htwist}[l(1)d(2)] 
!{\xcaph[1]@(0)}[ld] 
!{\xcaph[1]@(0)}[u(3)] 
!{\xcaph[1]@(0)}[ld] 
!{\htwist}[l(1)d(2)] 
!{\xcaph[1]@(0)}[u(3)] 
!{\xcaph[1]@(0)}[ld] 
!{\htwist}[l(1)d(2)] 
!{\xcaph[1]@(0)}[u(3)] 
} 
\else \vskip 5cm \fi
\caption{The braid $[[A_{1,2},A_{2,3}],[A_{2,3},A_{3,4}]]$ in $P_4$.}
\label{fig:brunnian}
\end{center}
\end{figure}

\begin{remark} 
When viewed as an element of the $4$-strand pure braid group, the braid $[[A_{1,2},A_{2,3}],[A_{2,3},A_{3,4}]] \in P_4$ is an example of a Brunnian braid.  The deletion of any strand trivializes the braid, see Figure \ref{fig:brunnian}. (Compare \cite[\S 3]{CFR10}.)
\end{remark}

\begin{remark}\label{rem:paris}
I.~Marin showed us an argument he credited to L.~Paris, showing the $P_5$ is not residually free, implying $P_n$ is not residually free for $n\geq 5$. Paris' argument uses the solution of the Tits conjecture for $B_5$ due to Droms, Lewin, and Servatius \cite{DLS90} (see also Collins \cite{Col94}) to produce a subgroup of $P_5$ isomorphic to the free product $\Z\ast (\Z \times F_2)$, as explained in \cite[Proposition 1.1]{Marin}. This latter group is not residually free, see \cite[Theorems 6 and~3]{Baum67}.
\end{remark}

Let $\varSigma$ be an orientable surface, possibly with punctures.  
Let $\varSigma^{\times n}=\varSigma \times \dots \times \varSigma$ denote the $n$-fold Cartesian product. 
The pure braid group $P_n(\varSigma)$ of the surface $\varSigma$ is the fundamental group of the configuration space
\[
F(\varSigma,n) = \{(x_1,\dots,x_n) \in \varSigma^{\times n} \mid x_i \neq x_j\ \text{if}\ i\neq j\}.
\]
of $n$ distinct ordered points in $\varSigma$.

\begin{corollary}
For $n \ge 4$, the pure braid group $P_n(\varSigma)$ is not residually free.
\end{corollary}
\begin{proof}
If $\varSigma \neq S^2$, the pure braid group $P_4$ embeds in $P_n(\varSigma)$, see Paris and Rolfsen \cite{PR99}.  
If $\varSigma=S^2$, then $P_4< P_n(S^2)$ for $n\ge 5$, and $P_4(S^2) \cong P_4/Z(P_4)$ (see, for instance, \cite{Bir75}). 
So the result follows from Theorem \ref{thm:notresfree}.
\end{proof}

\begin{remark} The fundamental group of the orbit space $F(\varSigma,n)/\Sigma_n$, where $\Sigma_n$ denotes the symmetric group, is the (full) 
braid group $B_n(\varSigma)$ of the surface $\varSigma$.  Recall from the Introduction that $B_3$ is not residually free, and that the only two-generator residually free groups are $\Z$, $\Z^2$, and $F_2$.  
If $\varSigma\neq S^2$ and $n\ge 3$, then $B_n(\varSigma)$ has a $B_3$ subgroup, so is not residually free.  
If $\varSigma=S^2$, then $B_3 < B_n(S^2)$ for $n\ge 4$, and $B_3(S^2) \cong B_3/Z(B_3)$. So $B_n(S^2)$ is not residually free for $n\ge 3$.  
\end{remark}

A complex hyperplane arrangement $\A=\{H_1,\dots,H_m\}$ is a finite collection of codimension one subspaces of $\C^n$.  Fix coordinates $(z_1,\dots,z_n)$ on $\C^n$, and for $1\le i \le m$, let $\ell_i(z_1,\dots,z_n)$ be a linear form with $\ker(\ell_i)=H_i$.  
The product $Q=Q(\A)=\prod_{i=1}^m \ell_i$ is a defining polynomial for $\A$.  
The group $G(\A)$ of the arrangement is the fundamental group of the complement 
$M(\A)=\C^n \smallsetminus \bigcup_{i=1}^m H_i= \C^n \smallsetminus Q^{-1}(0)$.

The arrangement $\A_{r,1,n}$ with defining polynomial
\[
Q=Q(\A_{r,1,n}) = z_1\cdots z_n \prod_{1\le i<j \le n}(z_i^r - z_j^r)
\]
is known as the full monomial arrangement (it is the reflection arrangement corresponding to the full monomial group $G(r,1,n)$).  Note that the arrangement $\A_{2,1,n}$ is the Coxeter arrangement of type $B_n$.  The complement $M(\A_{r,1,n})$ of the full monomial arrangement may be realized as the orbit configuration space
\[
F_{\Gamma}(\C^*,n) = \{(x_1,\dots,x_n) \in (\C^*)^{\times n} \mid \Gamma\cdot x_i \cap \Gamma\cdot x_j
=\emptyset\ \text{if}\ i\neq j\}
\]
of ordered $n$-tuples of points in $\C^*$ which lie in distinct orbits of the free action of $\Gamma=\Z_r$ on $\C^*$ by multiplication by the primitive $r$-th root of unity $\exp(2\pi\sqrt{-1}/r)$.

Call the fundamental group $P(r,1,n)=G(\A_{r,1,n})$ the pure monomial braid group.  For $n=1$, $P(r,1,1)\cong \Z$, and for $n=2$, it is well known that $P(r,1,2) \cong \Z\times F_{r+1}$.  Hence, $P(r,1,n)$ is residually free for $n\le 2$.

\begin{corollary}
For $n\ge 3$, the pure monomial braid group $P(r,1,n)$ is not residually free.
\end{corollary}
\begin{proof}
For $n\ge 3$, it follows from \cite{Coh01} that the pure braid group $P_4$ embeds in $P(r,1,n)$. 
So the result follows from Theorem \ref{thm:notresfree}.
\end{proof}

\begin{remark} The fundamental group of the orbit space $M(\A_{r,1,n})/G(r,1,n)$ is the (full) monomial braid group $B(r,1,n)$.  This group admits a presentation with generators 
$\rho_0,\rho_1,\dots,\rho_{n-1}$ and relations 
\[
(\rho_0\rho_1)^2=(\rho_1\rho_0)^2\!,\   
\rho_i\rho_{i+1}\rho_i=\rho_{i+1}\rho_i\rho_{i+1}\, (1\le 
i<n),\ \rho_i\rho_j=\rho_j\rho_i\  (|j-i|>1).
\]
Observe that $B(r,1,n)$ is independent of $r$, and is the Artin group of type $B_n$. 
For $n \ge 3$, the group $B(r,1,n)$ has a $B_3$ subgroup, so is not residually free.  The group $B(r,1,2)$ is not residually free, since it is a two-generator group which is not free or free abelian.
\end{remark}

Let $\varGamma$ be a Coxeter graph, with associated Artin group $A$ and pure Artin group $P$.  
We say that $\varGamma$ contains an $A_k$ subgraph if it contains a path of length $k$ with unlabelled edges as a vertex-induced subgraph.

\begin{corollary}
If $\varGamma$ contains an $A_3$ subgraph, then the associated pure Artin group $P$ is not residually free.
\end{corollary}
\begin{proof}
If $\varGamma$ contains an $A_3$ subgraph, it follows from van der Lek \cite{Lek} (see also \cite{Par97}) that $P$ has a $P_4$ subgroup. 
So the result follows from Theorem \ref{thm:notresfree}.
\end{proof}

Note that this includes the pure Artin groups associated to all irreducible Artin groups of finite type and rank at least $3$, except types $H_3$ and $F_4$.

\begin{remark}
If $\varGamma$ contains an $A_2$ subgraph, then the (full) Artin group $A$ has a $B_3$ subgroup, so is not residually free.  This includes all irreducible Artin groups of finite type and rank at least $2$, except type $I_2(m)$.
\end{remark}

\begin{ack} Portions of this project were completed during during the intensive research period ``Configuration Spaces: Geometry, Combinatorics and Topology,'' May-June, 2010, at the Centro di Ricerca Matematica Ennio De Giorgi in Pisa.  The authors thank the institute and the organizers of the session for financial support and hospitality, and for the excellent working environment. The authors are grateful to Luis Paris and Ivan Marin for helpful conversations.
\end{ack}

\newcommand{\arxiv}[1]{{\texttt{\href{http://arxiv.org/abs/#1}{{arXiv:#1}}}}}

\newcommand{\MRh}[1]{\href{http://www.ams.org/mathscinet-getitem?mr=#1}{MR#1}}

\newcommand{\vs}{\vspace{0pt}}

\bibliographystyle{amsalpha}

\end{document}